\documentclass[12pt]{article}
\usepackage{amsmath,amsfonts,amssymb,amsthm}
\usepackage[margin=1in]{geometry}
\usepackage[round]{natbib}
\theoremstyle{plain}
\newtheorem{theorem}{Theorem}[section]
\newtheorem{corollary}[theorem]{Corollary}

\theoremstyle{definition}
\newtheorem{assumption}[theorem]{Assumption}

\newtheorem{example}[theorem]{Example}
\theoremstyle{remark}

\numberwithin{equation}{section}

\begin{document}

\title{Stochastic Differential Equations on the Simplex}


\author{Paolo Guasoni\thanks{Dublin City University, School of Mathematical Sciences, D09 W6Y4 Dublin, Ireland, and Universit\`a di Bologna, Dipartimento di Statistica, Via Belle Arti 41, 40126 Bologna, Italy. \texttt{paolo@guasoni.it}. Partially supported by SFI (16/IA/4443, 16/SPP/3347).}\and Eberhard Mayerhofer\thanks{University of Limerick, Department of Mathematics and Statistics, V94 T9PX Limerick, Ireland. \texttt{eberhard.mayerhofer@ul.ie}}}

\maketitle

\begin{abstract}
We propose a framework for simplex-valued diffusion processes, formulated as solutions to stochastic differential equations (SDEs) driven by antisymmetric Brownian motion. Our edge-wise representation includes the multivariate Jacobi processes and the general class of continuous simplex-valued polynomial diffusions, and also allows bounded, locally Lipschitz state-dependent drift rates and edge weights. We establish weak solutions, and also global strong existence and pathwise uniqueness for strictly positive initial states under a sufficient boundary non-attainment condition.
\end{abstract}

\section{Introduction}

This paper presents a general class of stochastic differential equations whose solutions lie in the unit simplex.

Any diffusion $X$ whose state space is the unit interval yields a process $Z$ on the simplex $\Delta^2:=\{(x,y)\in\mathbb R^2\mid x,y\geq 0, x+y=1\}$, by setting $Z:=(X,1-X)$. One of the best-known and most tractable models of this type is the Wright--Fisher diffusion of \citet{kimura1954stochastic}, which arises as a large-population diffusion limit of the discrete Wright--Fisher model under an appropriate time scaling. The representation of such processes as stochastic differential equations became more widely known through models of stochastic interest rates and stochastic volatility (cf. \citealp{delbaen,veraart,ackerer}), in which $X$ is driven by a Brownian motion $W$,

\begin{equation}\label{eq: SDE1}
dX_t=(b-\beta X_t)dt+\sigma \sqrt{X_t(1-X_t)}dW_t,\quad X_0\in [0,1],
\end{equation}
and the parameters satisfy $0\leq b\leq \beta$ and $\sigma\geq 0$.

The polynomial character of this process has led researchers in finance, economics and other fields to generalize it to higher dimensions, obtaining tractable multifactor processes. First, higher-dimensional versions of the stochastic differential equation \eqref{eq: SDE1} were developed by \citet{gourieroux2006multivariate}. The multivariate Jacobi diffusion is a process \( X_t := (X^1_t, X^2_t, \ldots, X^m_t) \) satisfying the dynamics
\begin{equation}\label{eq: jacobm}
dX^k_t=-a(X^k_t-\pi^k)dt+\sqrt{c X_t^k}dW_t^k-X_t^k \sum_{l=1}^m\sqrt{c X_t^l} dW_t^l,
\end{equation}
with $a,c>0$, $\pi^j>0$ for all $1\leq j\leq m$, and $\sum_j \pi^j=1$. Furthermore, $W=(W^1,\dots,W^m)$ is an $m$-dimensional, standard Brownian motion. The resulting quadratic covariations are given by

\begin{equation}\label{eq: bill}
d\langle X^k,X^l\rangle_t=\begin{cases} c X_t^k(1-X_t^k) dt,\quad\quad k=l\\-c X_t^kX_t^l dt,\quad k\neq l\end{cases}.
\end{equation}

Second, \citet{filipovic2016polynomial} characterised the general class of simplex-valued polynomial diffusions, whose quadratic covariances are of the form
\begin{equation}\label{eq: fil}
d\langle X^k,X^l\rangle_t=\begin{cases} X^k_t\sum_{r=1}^m \Pi_{kr} X_t^r dt,\quad\quad k=l\\-\Pi_{kl} X_t^kX_t^l dt,\quad k\neq l\end{cases},
\end{equation}
where $\Pi_{kl}=\Pi_{lk}\geq 0$, with $\Pi_{kk}=0$ for $1\leq k,l\leq m$. Hence, the multivariate Jacobi diffusions \eqref{eq: jacobm} are a special case of polynomial ones, with a restricted parameter set (in particular, $\Pi_{ij}=c$ for all $1\leq i\neq j\leq m$).

Antisymmetric edge-noise representations have already appeared in the literature on interacting particle systems on graphs. In particular, \citet{itoh1998explicit} study a diffusion limit in which independent Brownian motions are associated with the edges of a graph and act with opposite signs at their endpoints. This corresponds to a drift-free special case of our formulation with constant edge weights that vanish away from the graph edges. \citet{mano2022measure} investigates the resulting measure-on-graph-valued diffusion, its stationary configurations and moment properties, and also considers a particular linear inward drift. Other related conservative noise constructions arise in finite-state mean-field games \citep{bayraktar2021finite} and in stochastic energy-exchange models \citep{peletier2014large}. In contrast, our framework permits bounded, locally Lipschitz state-dependent drift rates and edge weights, and our main result gives a sufficient condition for global strong solutions in the relative interior of the simplex. For constant coefficients, the representation covers the continuous polynomial simplex diffusions characterised by \citet{filipovic2016polynomial}.

\subsection*{Organisation of the paper}

In Section 2, we introduce a general class of stochastic differential equations on the simplex, driven by antisymmetric Brownian motion. Our main result establishes the existence of unique global strong solutions in the relative interior of the simplex under a sufficient condition for boundary non-attainment. In Section 3, we apply this framework to polynomial diffusions, establish the existence and uniqueness in law of weak solutions on the closed simplex under weaker conditions, and discuss pathwise uniqueness in special cases.

\section{SDE dynamics on the simplex}\label{sec: model}

Fix $m\geq 2$. Denote by \( \Delta^m \) the unit simplex in $\mathbb R^m$ (of affine dimension $m-1$), defined by

\[
\Delta^m = \left\{ (x_i)_{i=1}^m \in \mathbb{R}^m \mid x_i \geq 0 \; \text{for all} \; i, \; \sum_{i=1}^m x_i = 1 \right\}
\]
and its relative interior
\[
\Delta^m_\circ = \left\{ (x_i)_{i=1}^m \in \mathbb{R}^m \mid x_i > 0 \; \text{for all} \; i, \; \sum_{i=1}^m x_i = 1 \right\}.
\]

To describe the stochastic evolution of simplex-valued processes, consider a probability space $(\Omega,\mathcal F,\mathbb P)$, supporting $m(m-1)/2$ independent standard Brownian motions $Z^{ij}$, $1\leq i<j\leq m$. Write $Z_t=(Z_t^{ij})_{1\leq i<j\leq m}$. Denote by $(\mathcal F_t)_{t\geq 0}$ the augmentation of their natural filtration.

The Brownian motions $Z$ are used to construct the antisymmetric $m\times m$ matrix Brownian motion\footnote{
The antisymmetric component $(B-B^\top)/\sqrt{2}$ of a matrix $B$ of independent standard Brownian motions has standard Brownian entries above the diagonal.

Antisymmetric Brownian motion has also been used for matrix-valued processes and processes on ellipsoids \citep{norris1986brownian}; for earlier simplex-valued models driven by antisymmetric edge noise, see \citet{itoh1998explicit} and \citet{mano2022measure}.} $W$ defined by $W^{ij}=Z^{ij}$ and $W^{ji}=-W^{ij}$ for $1\leq i<j\leq m$, with vanishing diagonal processes $W^{ii}=0$ for $1\leq i\leq m$.

We consider solutions $U$ of the $m$-dimensional stochastic differential equation

\begin{equation}\label{eq: SDE for simplex under P}
dU_t^i =\sum_{j=1}^m U_t^j\beta^i_j (U_t) dt+\sum_{j=1}^m\sqrt{ U^i_t U^j_t\Pi_{ij}(U_t)}dW_t^{ij}, \quad U^i_0=u^i,\quad 1\leq i\leq m,
\end{equation}
where the drift and diffusion coefficients satisfy the following conditions.

\begin{assumption}\label{ass: 1 gen}

For $1\le i,j\le m$, let $\Pi_{ij}$ and $\beta^i_j$ be bounded, locally Lipschitz functions from $\Delta^m_\circ$ to $\mathbb R$, satisfying

\begin{enumerate}

\item \label{ass1}

$\sum_{k=1}^m\beta^k_i=\Pi_{ii}=0$ for all $1\leq i\leq m$ and $\Pi_{ij}=\Pi_{ji}$ for all $1\leq i\ne j\leq m$;

\item \label{ass3} For every $1\leq i<j\leq m$, either $\Pi_{ij}$ vanishes identically or $\Pi_{ij}(u)>0$ for every $u\in\Delta^m_\circ$.

\end{enumerate}

\end{assumption}

The main result establishes the existence of global solutions in the relative interior of the unit simplex under a mild condition on the drift.

\begin{theorem}\label{th: main}

Suppose

\begin{equation}\label{eq: boundary non-attainment}
\beta^i_j(u)\geq \Pi_{ij}(u)/2, \quad u\in \Delta^m_{\circ},\quad 1\leq i\ne j\leq m.
\end{equation}

Then a unique global strong solution of \eqref{eq: SDE for simplex under P} exists on $ \Delta^m_{\circ}$. Moreover, for any initial value $u\in \Delta^m_\circ$, the relative boundary of the simplex $\Delta^m$ is not attained in finite time, almost surely.

\end{theorem}
\begin{proof}
Fix $u\in\Delta^m_\circ$ and use $m-1$ affine coordinates on the hyperplane $\sum_i x^i=1$. Since every nonzero $\Pi_{ij}$ is strictly positive in the relative interior, the diffusion coefficients $\sqrt{x^ix^j\Pi_{ij}(x)}$ are locally Lipschitz there. The remaining coefficients are also locally Lipschitz. Thus the equation admits a pathwise unique strong solution up to its maximal lifetime in $\Delta^m_\circ$. The antisymmetry of $W$ and the column-sum identity for $\beta$ preserve the total mass.

Set $Z_t=\prod_{i=1}^m U_t^i$ and $h(z)=\log z$ for $z>0$.
For $t$ before the first exit from $\Delta^m_\circ$, It\^o's formula gives
\begin{equation}\label{eq: log-drift}
 dh(Z_t)=\left[\sum_{i=1}^m\beta^i_i(U_t)
 +\sum_{i=1}^m\sum_{j\ne i}\frac{U_t^j}{U_t^i}
 \left(\beta^i_j(U_t)-\frac{\Pi_{ij}(U_t)}2\right)\right]dt+d\mathcal M_t,
\end{equation}
where $\mathcal M$ is a continuous local martingale on the stochastic interval before exit. By \eqref{eq: boundary non-attainment}, the finite-variation part of \eqref{eq: log-drift} is bounded from below on every finite time interval: indeed, its drift is at least
\[
 \sum_{i=1}^m\beta^i_i(U_t)\geq -C,\qquad
 C:=\sum_{i=1}^m\|\beta^i_i\|_\infty<\infty.
\]
If the solution reached the relative boundary in finite time, then by continuity and conservation of total mass $Z_t\downarrow0$ along the approach to the exit time, and hence $h(Z_t)\to-\infty$. McKean's argument, in the form of \citet[Proposition~4.3]{PMS}, rules this out: a continuous local martingale cannot converge to $-\infty$ at a finite endpoint while the finite-variation part is bounded from below. Thus the relative boundary is not attained in finite time. Since the simplex is bounded and the coefficients are locally Lipschitz in its relative interior, there is no other finite-time explosion. The solution is therefore global and pathwise unique.
\end{proof}

\begin{example}\label{ex: power coeff}
Let $\alpha,\beta,\gamma\geq0$, and let $b^i_j$ and $p_{ij}$ be nonnegative constants for $i\ne j$, with $p_{ij}=p_{ji}$. Define
\[
\Pi_{ij}(u)=p_{ij}(u^i)^\gamma(u^j)^\gamma,
\qquad
\beta^i_j(u)=b^i_j(u^i)^\alpha(u^j)^\beta,
\qquad i\ne j,
\]
and set $\Pi_{ii}(u)=0$ and
\[
\beta^i_i(u)=-\sum_{k\ne i}\beta^k_i(u),\qquad 1\leq i\leq m.
\]
Then Assumption \ref{ass: 1 gen} holds. The non-attainment condition \eqref{eq: boundary non-attainment} is satisfied if
\[
\max(\alpha,\beta)\leq\gamma,
\qquad b^i_j\geq p_{ij}/2,\qquad i\ne j.
\]
For example, if $\gamma=1$, $\alpha=\beta=0$ and $p_{ij}=c_{ij}^2$ with $c_{ij}\geq0$, the dynamics become
\[
dU_t^i=\sum_{j=1}^mU_t^j\beta^i_j\,dt
+\sum_{j=1}^m c_{ij}U_t^iU_t^j\,dW_t^{ij},
\qquad U_0^i=u^i,\quad 1\leq i\leq m.
\]
The case $\alpha=\beta=\gamma=0$ is discussed in the next section.
\end{example}

\section{Polynomial processes}\label{sec: Polynomial}

For constant coefficients $\Pi$ and $\beta$ (with linear and quadratic terms factored out in \eqref{eq: SDE for simplex under P}), the SDE provides a representation of the polynomial diffusions of \citet{filipovic2016polynomial}. In this case, Assumption \ref{ass: 1 gen}(ii) is always satisfied, as either the constant coefficient $\Pi_{ij}$ vanishes, or is strictly positive. This special case also corresponds to Example \ref{ex: power coeff} for $\alpha=\beta=\gamma=0$.

\subsection{Boundary non-attainment}

The SDE \eqref{eq: SDE for simplex under P} takes the simpler form

\begin{equation}\label{eq: SDE constant coefficients}
dU_t^i =\sum_{j=1}^m U_t^j\beta^i_j dt+\sum_{j=1}^m\sqrt{ U^i_t U^j_t\Pi_{ij}}dW_t^{ij}, \quad U^i_0=u^i,\quad 1\leq i\leq m
\end{equation}

The parameterisation in \citet{filipovic2016polynomial} also includes a constant drift $b=(b^i)_{1\leq i\leq m}$. Since $\sum_{j=1}^m u^j=1$ on the simplex, this term can be absorbed into the linear drift by setting $\beta^i_j=b^i+\widetilde\beta^i_j$, because
\begin{equation}\label{eq: comp}
b^i+\sum_{j=1}^m\widetilde\beta^i_j u^j
=\sum_{j=1}^m(b^i+\widetilde\beta^i_j)u^j
=\sum_{j=1}^m\beta^i_j u^j.
\end{equation}
We use the parameterisation \eqref{eq: SDE constant coefficients} throughout.

Theorem \ref{th: main} now gives the following corollary directly. Its equality cases follow from our proof; a direct invocation of the boundary criteria in \citet[Theorem 5.7]{filipovic2016polynomial} would require checking the hypotheses of that theorem separately.

\begin{corollary}\label{cor: strong}

For every initial value $u\in\Delta^m_\circ$, the SDE \eqref{eq: SDE constant coefficients} has a unique global strong solution if

\begin{equation}\label{eq: nonattainment}
\beta^i_j\geq \Pi_{ij}/2\quad \text{for all}\quad 1\leq i\neq j\leq m.
\end{equation}

The solution does not hit the boundary $\Delta^m\setminus \Delta^m_\circ$ in finite time, almost surely.

\end{corollary}

\subsection{Weak solutions}

For the constant-coefficient system \eqref{eq: SDE constant coefficients}, assume throughout this subsection that $\Pi_{ij}=\Pi_{ji}\geq0$, $\Pi_{ii}=0$, $\sum_i\beta^i_j=0$ for each $j$, and
\begin{equation}\label{eq: inward}
 \beta^i_j\geq0,\qquad i\ne j.
\end{equation}
The nonnegative off-diagonal drift condition is also necessary for a simplex-valued solution from \emph{every} starting point of the closed simplex. Indeed, starting from the vertex $e_j$ and taking expectations gives
\[
 0\leq\mathbb E_{e_j}[U_t^i]=t\beta^i_j+o(t),\qquad i\ne j,
\]
so $\beta^i_j\geq0$. The strict condition $\beta^i_j>0$ used in \citet[Theorem 5.3]{filipovic2016polynomial} ensures additional boundary-occupation properties, but strict positivity is not needed for weak existence; see also the Type~0 case in \citet[Section 6]{cuchiero2018jump}.

To prove sufficiency, write $x_i^+=\max(x_i,0)$ and extend the drift and edge-noise coefficients to all of $\mathbb R^m$ by
\begin{equation}\label{eq: ambient-extension}
 \widehat b_i(x)=\sum_{j=1}^m\beta^i_jx_j^+,
 \qquad
 \widehat g_{ij}(x)=\sqrt{\Pi_{ij}x_i^+x_j^+}\,(e_i-e_j),\quad i<j.
\end{equation}
With independent standard Brownian motions $W^{ij}$ for $i<j$, consider
\[
 dU_t=\widehat b(U_t)\,dt+\sum_{i<j}\widehat g_{ij}(U_t)\,dW_t^{ij},
 \qquad U_0=u\in\Delta^m.
\]
These coefficients are continuous and have at most linear growth, hence a global weak solution exists on a suitable stochastic basis (cf. \citealp{stroock1997multidimensional,protter}). If $x_i\leq0$, every noise coefficient in component $i$ vanishes and
\[
 \widehat b_i(x)=\sum_{j\ne i}\beta^i_jx_j^+\geq0.
\]
Thus the half-space invariance criterion \citep[Lemma B.1(ii)]{fm09} applies to each coordinate half-space, and $U_t^i\geq0$ for all $t\geq0$ almost surely. Each edge-noise vector has zero coordinate sum, and
$\sum_i\widehat b_i(x)=\sum_j(\sum_i\beta^i_j)x_j^+=0$. Consequently $\sum_iU_t^i=1$ and $U_t\in\Delta^m$ for all $t\geq0$. On this simplex, \eqref{eq: ambient-extension} agrees with \eqref{eq: SDE constant coefficients}.

The resulting simplex-valued law is unique: its generator preserves the finite-dimensional spaces of polynomials of degree at most $n$, for every $n$. The conditional polynomial moments therefore satisfy finite-dimensional linear equations; polynomial density on the compact simplex determines the transition kernels and the path law (cf. \citealp[Theorem 4.2 and Remark 6.7]{filipovic2016polynomial}). Under weak inequalities in \eqref{eq: inward}, boundary faces may be absorbing, so zero occupation time at the boundary is not asserted.

\subsection{Pathwise uniqueness}

Corollary \ref{cor: strong} gives strong existence and pathwise uniqueness for interior initial states under the non-attainment condition. For a solution that starts on, or may reach, the relative boundary, further arguments are needed. We discuss the one-dimensional case and the case of vanishing drift, and distinguish a diagonal system from our edge-driven equation.

For $m=2$, setting $X_t=U_t^1$ and $U_t^2=1-X_t$ gives
\[
 dX_t=\bigl[\beta^1_2+(\beta^1_1-\beta^1_2)X_t\bigr]dt
 +\sqrt{\Pi_{12}X_t(1-X_t)}\,dW_t^{12}.
\]
Under the inward drift condition \eqref{eq: inward}, the equation has weak solutions on $[0,1]$, and its diffusion coefficient is H\"older continuous of order $1/2$. The criterion of \citet[Theorem 1]{yamada1971uniqueness} yields pathwise uniqueness, and hence strong existence, also for boundary initial states.

For comparison, consider the different system with diagonal diffusion coefficients
\[
 d\widetilde U_t^i
 =\sum_{j=1}^m\beta^i_j\widetilde U_t^j\,dt
 +\sigma\sqrt{\max\{\widetilde U_t^i(1-\widetilde U_t^i),0\}}\,dB_t^i,
\]
where $B^1,\ldots,B^m$ are independent standard Brownian motions. The coordinatewise H\"older criterion applies to this diagonal system, giving pathwise uniqueness. It is \emph{not} a representation of \eqref{eq: SDE constant coefficients}: its cross-variations vanish, whereas the edge-driven simplex equation has
$d\langle U^i,U^j\rangle_t=-\Pi_{ij}U_t^iU_t^j\,dt$ for $i\ne j$. In particular, their laws do not generally agree.

Finally, suppose $\beta^i_j=0$ for all $i,j$ and retain arbitrary constant symmetric nonnegative edge weights. The edge-driven SDE has a unique strong solution from every $u\in\Delta^m$, and pathwise uniqueness holds for solutions driven by the same edge Brownian motions. To see this, start on the relative interior of the face supported by $S=\{i:u^i>0\}$. The coefficients corresponding to the positive coordinates are locally Lipschitz there; all other coordinates are set to zero. Solve strongly up to the first time one or more positive coordinates hit zero. On this interval all coordinates are bounded martingales, and the bounded noise coefficients ensure that the solution has a continuous limit at a finite exit time. A coordinate that reaches zero cannot later become positive in any simplex-valued solution: it is a nonnegative martingale, so its conditional expectation after the hitting time remains zero. Hence the process continues on the smaller face, using the same Brownian increments. Each finite face change removes at least one active coordinate; there can be at most $m-1$ such changes. If the lifetime on a face is infinite, the construction already gives the required continuation. Local pathwise uniqueness on each face and induction over face changes establish pathwise uniqueness for the full equation.

This does \emph{not} imply eventual absorption at a vertex. For example, let $m=3$, $\Pi_{12}=\Pi_{23}=1$ and $\Pi_{13}=0$. Every point $(x_1,0,x_3)$ with $x_1+x_3=1$ is stationary, and the process may fail to hit any vertex. Such stationary configurations are related to independent sets of the interaction graph; see \citet{itoh1998explicit} and \citet{mano2022measure}. The preceding face-by-face argument does not resolve pathwise uniqueness for arbitrary nonzero inward drift when \eqref{eq: nonattainment} fails.
%
%


\end{document}